\documentclass{amsart}
\usepackage{amsfonts}

\theoremstyle{plain}

\newtheorem{corollary}{Corollary}

\newtheorem{lemma}{Lemma}

\newtheorem{proposition}{Proposition}

\numberwithin{equation}{section}

\begin{document}
\title[Invariance principles]{Invariance principles for random sums of random variables}
\author{Gane Samb LO}


\begin{abstract}
This note investigates invariance principles for sums of N(nt) iid radom variables, where n is an integer, t is a positive real number and N(u) is a stochastic process with nonnegative integer values. We show that the sequence of sums of these random variables denoted S(n,t), when appropriately centered and normalized, weakly converges to a Gaussian process. We give sufficient conditions depending on the expectation of N(nt) which allows to rescale S(n,t) into a stochastic S(n,a(t)) weakly converging to a Brownian motion.
\end{abstract}
\maketitle
\Large

\section{Introduction}

\label{sec1}

Let $X_{1},X_{2},...$ be sequence of real random variables defined on the same probability space $(\Omega ,\mathcal{A},\mathbb{P})$ and suppose that
this probability space holds a stochastic process $(N(t))_{t\geq 0}$ taking its values in the set $\mathbb{N}$  of nonnegative integers with $N(0)=0$. Define for $t\geq 0$, \begin{equation*}
S(t)=\sum_{h=1}^{N(t)}X_{h},\text{ for }N(t)\geq 1,
\end{equation*}

\noindent and $S(t)=0$ for $N(t)<1$ and for any $n\geq 1$ and $t\geq 0,$%
\begin{equation*}
S_{n}(t)=S(nt).
\end{equation*}

\noindent This note investigates possible invariance principles for the
normalized random sums $\{S_{n}(t)/\sqrt{c_{n}},t\geq 0\}$ , where $\left(
c_{n}\right) _{n\geq 1}$ is a sequence of positive constants to be defined
later, in the space of locally bounded functions $\ell ^{\infty }$, that is
the space functions which are bounded on compact sets $[0,T]$. Indeed for
each $T>0$, the sequence

\begin{equation*}
\{S_{n}(t)/\sqrt{c_{n}},0\leq t\leq T\}=\{\sum_{h=1}^{N(nt)}X_{h}/\sqrt{c_{n}%
},0\leq t\leq T\}
\end{equation*}

\noindent is in the space $\ell ^{\infty }(T)=\ell ^{\infty }([0,T])$ of
real bounded functions on $[0,T]$ equipped with the sup-norm%
\begin{equation*}
\left\Vert x\right\Vert _{\infty ,T}=\sup \sup_{t\in \lbrack 0,T]}\left\vert
x(t)\right\vert ,x\in \ell ^{\infty }(T),
\end{equation*}

\noindent since for any $\omega \in \Omega $ for any $n\geq 1$%
\begin{equation*}
\left\Vert S_{n}\right\Vert _{\infty ,T}=\sup_{t\in \lbrack 0,T]}\left\vert
S_{n}(t)\right\vert \leq \max_{1\leq h\leq N(T)(\omega )}\left\vert
\sum_{h=1}^{k}X_{h}(\omega )/\sqrt{c_{n}}\right\vert <+\infty .
\end{equation*}

\noindent So it will be enough to study the weak convergence of the
stochastic processes $\{S_{n}(t),0\leq t\leq T\},$ for each fixed $T>0$.%
\newline

\noindent We intend to proceed to a general study from the case of
independent variables $X_{1},X_{2},...$ to dependent data. The case of
associated sequences $X_{1},X_{2},...$ would place this study in an
currently active research field.\newline

But for the beginning, we explore the case where the $X_{1},X_{2},..$ are $%
iid$ centered random variables with finire second moments (taken to be one).%
\newline

\noindent Let us make some hypotheses we may need.\newline

\noindent \textbf{(HX1A)} The random variables $X_{1},X_{2},..$ are centered
iid with finite variance ($\sigma ^{2}=1)$ with common characteristic
function $\varphi$.\\

\noindent \textbf{(HN1A)} The stochastic process $(N(t))_{t\geq 0}$ is
everywhere increasing and has independent increments such that for any $%
0\leq s<t,$ 
\begin{equation*}
0\leq \mathbb{E}(N(nt)-N(ns))=a_{n}(s,t)\in R\text{ and }a_{n}(s,t)%
\rightarrow +\infty \text{ as }n\rightarrow +\infty .
\end{equation*}

\noindent \textbf{(HN2A)} For each $0\leq s<t,$ the moment function of $%
\Delta N_{n}(s,t)$ is well defined on $\mathbb{R}$ and the law of large
numbers 
\begin{equation*}
\frac{\Delta N_{n}(s,t)}{a_{n}(s,t)}\rightarrow _{P}1,
\end{equation*}

\noindent holds.\newline

\noindent \textbf{(HN3A)} There exists a sequence ($c_{n}>0)_{n\geq 0}$ such
that for any $0\leq s<t,$ $a_{n}(s,t)/c_{n}$ converges to a positive real
number $\Delta a(s,t).$.\newline

\noindent \textbf{(HTA)} For some $p\geq 1,$ for any $t\in \lbrack 0,T],$%
\begin{equation*}
\lim_{\delta \rightarrow 0}\frac{1}{\delta }\lim \sup_{n\rightarrow +\infty }%
{\large \mathbb{E}}\left( \left\{ \frac{{\large N(nt+n\delta )-N(nt)}}{c_{n}}%
\right\} ^{p/2}\right) {\large =0}
\end{equation*}

\noindent and the sequence 
\begin{equation*}
\mathbb{E}\left( \left\vert \frac{S_{r}}{\sqrt{r}}\right\vert ^{p}\right)
,p\geq 1,
\end{equation*}

\noindent is bounded.\\

\noindent \textbf{(HA1)} The fonction $a(\circ )$ is strictly increasing and invertible, and $%
a^{-1}(\circ )$ is a Lipschitz function withe index $k>0.$

\section{Main results}

\label{sec2}

The main results in the iid case are the following.

\begin{proposition}
\label{propA} Assume (HX1A),(HX2A), (HN1A), (HN2A), (HN3A) ,
(HTA)  hold. Then, for each $T>0,$ the sequence of stochastic processes $%
S_{n,T}=\{S_{n}(t)/\sqrt{c_{n}},0\leq t\leq T\}$ weakly converges, in the
space $\ell ^{\infty }([0,T])$ to a centered Gaussian process 
$$
\mathbb{G}_{T}=\{\mathbb{G},0\leq t\leq T\}
$$ 

\noindent with covariance function defined by

\begin{equation*}
\Gamma _{\mathbb{G}}(s,t)=a(\min (s,t)),(s,t)\in \lbrack 0,T]^{2},
\end{equation*}

\noindent where 
\begin{equation*}
a(u)=\Delta a(0,u),u\geq 0,
\end{equation*}

\noindent in the sense that for any continuous and bounded function $%
f:\left( \ell ^{\infty }(T),\left\Vert \circ \right\Vert _{\infty ,T}\right)
\rightarrow \mathbb{R},$%
\begin{equation}
\mathbb{E}f(S_{n,T})\rightarrow \mathbb{E}f(\mathbb{G}_{T})\text{ as }%
n\rightarrow +\infty .  \label{defConV}
\end{equation}
\end{proposition}

\bigskip
\bigskip \noindent The function $a(u),u\geq 0,$ is non-increasing. If $%
a(\circ)$ is increasing and invertible, we will be able to rescale the time
and get:

\begin{corollary}
\label{corA} Let $a(\circ )$ be an increasing and invertible function. If (HA1) holds in
addition of the assumptions of Proposition \ref{propA}, then the sequence of stochastic processes $%
\{S_{n}(a^{-1}(t)))/\sqrt{c_{n}},0\leq t\leq T\}$ weakly converges, in the
space $\ell ^{\infty }([0,T])$ to a centered Gaussian process 
$$
\{\mathbb{G(}a^{-1}(t)),0\leq t\leq T\}^{{}},0\leq t\leq T\},
$$

\noindent where $\mathbb{G(}%
a^{-1}(\circ ))=B(\circ )$ Brownian motion on $[0,T]$
\end{corollary}

\bigskip \noindent \textbf{FIRST REMARKS}.\\

\noindent \textbf{(1)} The second part of Assumption (HTA), that is the sequence $\mathbb{E}%
\left( \left\vert S_{r}/\sqrt{r}\right\vert ^{p}\right) $ is bounded for
some $p>2$, seems to too strong. Indeed it implies that the random variables 
$X_{i}$ have a $pth$ finite moment. It is satisfied if $\mathbb{E}X_{i}^{4}$
is finite. Actually, this condition comes from the use the submartingale
argument. Other arguments should be investigated to avoid it.\\

\noindent \textbf{(2)} All the conditions on $N$ are easily satisfied if  $N$ is the counting
process of a homogenuous Poisson Process with parameter $\lambda $ with $%
c_{n}=n$ and $a(t)=\lambda t$.\\

\bigskip \noindent \textbf{PERSPECTIVES}.\\

\noindent \textbf{(1)}  : Deepen these results and weakening the conditions.\\

\noindent \textbf{(2)}  : Explore possible applivations.\\

\noindent \textbf{(3)}  : Explore the Hungarian approximations of the form
$$
\sup_{t \in [0,T]} |S_{n,T}(t)-\mathbb{G}_T(t)| = O(a_n) \ \ a.s. 
$$.

\noindent \textbf{(4)}  : Generalize results for associated random variables.\\

\noindent \textbf{(5)}  : Generalize results for depedent data.\\

\section{Proofs}

\subsection{Proof of Proposition \protect\ref{propA}}

\noindent Here, we use the theory of weak convergence of applications with
values in $\left( \ell ^{\infty }(T),\left\Vert \circ \right\Vert _{\infty
,T}\right) $ now popularized by the book of var der Vaart and Wellner \cite%
{vaart}. This theory allows to avoid the Skorohod metric at the cost of
using outer and inner integrals or probabilities. In the the following, $%
\mathbb{P}^{\ast }(B)$ is the outer probability of any subset $B\subset
\Omega ,$ defined by%
\begin{equation*}
\mathbb{P}^{\ast }(B)=\inf \{\mathbb{P}(A),\text{ A measurable, }B\subset
A\}.
\end{equation*}

\noindent By Theorem 1.5.4 and 1.5.6 in \cite{vaart}, $S_{n,T}/\sqrt{c_{n}}$
weakly converges to $\mathbb{G}_{T}$ in $\left( \ell ^{\infty
}(T),\left\Vert \circ \right\Vert _{\infty ,T}\right) $ if and only if%
\newline

\noindent \textbf{(a)} the finite distributions of $\left\{ S_{n,T}/\sqrt{%
c_{n}}\right\} _{n\geq 1}$ weakly converge to those of $\mathbb{G}_{T}$,%
\newline

\noindent and\\

\noindent \textbf{(b)} the sequence $\left\{ S_{n,T}/\sqrt{c_{n}}\right\}
_{n\geq 1}$ is asymptoticaly tight.\newline

\noindent From Theorem 1.5.6 in \cite{vaart}, from the adaptation of \cite%
{loTensTool} for Theorem 8.3 \ \ in Billinsgley \cite{billingsley}, it comes that,
if the finite distributions already converge, then $\left\{ S_{n,T}\right\}
_{n\geq 1}$ is asymptotically tight whenever

\begin{equation}
\lim_{\delta \rightarrow 0}\sup_{s\in \lbrack 0,T]}\lim \sup_{n}\frac{1}{%
\delta }{\large \mathbb{P}}^{\ast }{\large (\sup_{s-\delta <t<s+\delta ,t\in
\lbrack 0,T]}\frac{\left\vert S_{n}(s)-S_{n}(t)\right\vert }{\sqrt{c_{n}}}%
>\eta )=0.}  \label{formulaTensIndep}
\end{equation}

\noindent Based on these remarks, we are going to proceed into three steps.%
\newline

\noindent \textbf{First}, we need this lemma. Next, we will prove the points
(a) and (b) above.

\begin{lemma}
\label{lem01} Let $0=t_{0}<t_{1}<...<t_{k}=T>0,$ and $k\geq 2.$ Denote $%
N_{n}(t_{j})=N(nt_{j}),$ $j=0,...,k$ and $\Delta
_{n}N(t_{j})=N_{n}(t_{j})-N_{n}(t_{j-1}),$ $j=1,...,k.$ Set%
\begin{equation*}
Y_{j,n}=\sum_{h=N(t_{j-1})+1}^{N(t_{j})}X_{h},\text{ }j=1,...,k.
\end{equation*}

\noindent Suppose that $(HX1A)$ and $(HN1A)$ and $(HN1B)$ hold. Then \ the
random variables $Y_{j,n},$ $j=1,...,k$, are independent.
\end{lemma}

\bigskip \noindent \textbf{Proof of Lemma \ref{lem01}}. Put $N^{\ast
}=(N_{n}(t_{1}),...,N_{n}(t_{k})).$ The characteristic function of $%
(Y_{1,n},...,Y_{k,n})$ is 
\begin{eqnarray*}
\psi _{(Y_{1,n},...,Y_{k,n})}(v_{1},...,v_{k}) &=&\mathbb{E}\exp
(i\sum_{j=1}^{k}v_{j}Y_{j,n}) \\
&=&\mathbb{E}\prod\limits_{j=1}^{k}\exp
(iv_{j}\sum_{h=N(t_{j-1})+1}^{N(t_{j})}X_{h}).
\end{eqnarray*}

\noindent But, for $n^{\ast }=(n_{1},...,n_{k})$ such that $n_{0}=0\leq
n_{1}\leq ...\leq n_{k},$ we have

\begin{eqnarray*}
\mathbb{E}\left( \left( \prod\limits_{j=1}^{k}\exp
(i\sum_{h=N(t_{j-1})+1}^{N(t_{j})}X_{h})\right) /N^{\ast }=n^{\ast }\right)
&=&\mathbb{E}\left( \prod\limits_{j=1}^{k}\exp
(iv_{j}\sum_{h=n_{j-1}+1}^{n_{j}}X_{h})\right) \\
&=&\prod\limits_{j=1}^{k}\varphi (v_{j})^{n_{j}-n_{j-1}},
\end{eqnarray*}

\noindent so that 
\begin{equation*}
\mathbb{E}\left( \left( \prod\limits_{j=1}^{k}\exp
(i\sum_{h=N(t_{j-1})+1}^{N(t_{j})}X_{h})\right) /N^{\ast }\right)
=\prod\limits_{j=1}^{k}\varphi (v_{j})^{\Delta _{n}N(t_{j})}.
\end{equation*}

\noindent We get 
\begin{eqnarray*}
\psi _{(Y_{1,n},...,Y_{k,n})}(v_{1},...,v_{k}) &=&\mathbb{E}\left( \mathbb{E}%
\left( \left\{ \exp (i\sum_{j=1}^{k}v_{j}Y_{j,n})\right\} /N^{\ast }\right)
\right) \\
&=&\mathbb{E}\left( \prod\limits_{j=1}^{k}\varphi (v_{j})^{\Delta
_{n}N(t_{j})}\right) \\
&=&.\mathbb{E}\left( \prod\limits_{j=1}^{k}\exp (\Delta _{n}N(t_{j})\log
\varphi (v_{j}))\right) \\
&=&\mathbb{E}\prod\limits_{j=1}^{k}\psi _{\Delta _{n}N(t_{j})}(\log \varphi
(v_{j})).
\end{eqnarray*}

\noindent So we have the independance between the $Y_{j,n},j=1,...,k.$

\bigskip \noindent \textbf{Secondly}, let us address the weak convergences
of the finite distributions of $\left\{ S_{n,T}\right\} _{n\geq 1}.$ We have

\begin{proposition}
\label{prop01} \bigskip Suppose that $(HX1A)$ and $(HN1A)$ and $(HN1B)$
hold. Let $0=t_{0}<t_{1}<...<t_{k}=T>0,$ and $k\geq 2.$ Denote $%
N_{n}(t_{j})=N(nt_{j}),$ $j=0,...,k$ and $\Delta
_{n}N(t_{j})=N_{n}(t_{j})-N_{n}(t_{j-1}),$ $a_{n}(t_{j})=E\Delta
_{n}N(t_{j}),$ $j=1,...,k.$ Set%
\begin{equation*}
Y_{j,n}=\sum_{h=N(t_{j-1})+1}^{N(t_{j})}X_{h},\text{ }j=1,...,k.
\end{equation*}

\noindent We have

\noindent \textbf{(1)} The $k$-tuple ($Y_{1,n}/\sqrt{a_{n}(t_{1})}%
,...,Y_{k,n}/\sqrt{a_{n}(t_{k})})$ weakly converges to a centered $k$%
-Gaussian random vector $\Delta B=(\Delta B(1),\Delta B(2),...,\Delta B(k))$
with independent standard gaussian components.\newline

\noindent \textbf{(2)} Suppose there exists a sequence of positive numbers $%
c_{n}$ such that for any $0\leq s<t,a_{n}(s,t)/c_{n}\rightarrow \Delta
a(s,t).$ Denote $a(t_{j})=a(t_{j-1},t_{j}),j=1,...,k$.\\

\noindent  Then the $k$-tuple $(S_{n}(t_{1})/\sqrt{c_{n}},...,(S_{n}(t_{2})/\sqrt{c_{n}})$ weakly converges
to the Gaussian random vector 
\begin{equation} \small
(\Delta B(1)\sqrt{a(t_{1})},...,\Delta B(1)\sqrt{a(t_{1})}+\Delta B(2)\sqrt{%
a(t_{2})}+...+\Delta B(k)\sqrt{a(t_{k})})  \label{fdl01}
\end{equation}

\noindent with covariance matrix $\sum =(\sigma _{uj})_{1\leq i,j\leq k}$
such that

\begin{equation*}
\sigma _{ij}=\sum\limits_{h=1}^{\min (i,j)}a(t_{h}).
\end{equation*}
\end{proposition}

\noindent \textbf{Proof of Proposition \ref{prop01}}. By Lemma \ref{lem01},
we have that the $Y_{j,n}$ are independent. So we only need to establish the
weak convergence of each component to get the joint weak convergence. We
have for each $j\in \{1,2,...,k\},$ 
\begin{equation*}
\mathbb{E}(Y_{j,n})=\mathbb{E}(\mathbb{E}(Y_{j,n}/N^{\ast }=n^{\ast
})=0\times (n_{j}-n_{j-1})=0
\end{equation*}

\noindent and 
\begin{eqnarray*}
\mathbb{E}(Y_{j,n}^{2}) &=&\mathbb{E}(\mathbb{E}(Y_{j,n}/N^{\ast }=n^{\ast })
\\
&=&\mathbb{E}(\left\{
\sum_{h=n_{j-1}+1}^{n_{j}}X_{h}^{2}+\sum_{n_{j-1}+1\leq h\neq \ell \leq
n_{j}}X_{h}X_{\ell }\right\} =(n_{j}-n_{j-1}).
\end{eqnarray*}

\bigskip \noindent Then then $Y_{j,n}$ are centered and have variance $%
\mathbb{E}\Delta _{n}N(t_{j})=a_{n}(t_{j}).$ Let us show that each $Y_{j,n}/%
\sqrt{a_{n}(t_{j})}$ converges to a $N(0,1)$ random variable. By the
computations in the Lemma \ref{lem01}, we have 
\begin{equation*}
\psi _{Y_{j,n}/\sqrt{a_{n}(t_{j})}}(v)=E\exp (ivY_{j,n}/\sqrt{a_{n}(t_{j})}%
)=E\exp (\Delta _{n}N(t_{j})\log \varphi (v/\sqrt{a_{n}(t_{j})})).
\end{equation*}

\noindent Let us use the uniform expansion of 
\begin{equation*}
\sup_{\left\vert v\right\vert \leq u}v^{-3}\left\vert \varphi (v)-1+\frac{%
v^{2}}{2}\right\vert =A(u)\text{ with }\lim \sup_{u\rightarrow 0}A(u)<+\infty
\end{equation*}

\noindent and 
\begin{equation*}
\sup_{\left\vert v\right\vert \leq u}v^{-2}\left\vert \log
(1+v)-v\right\vert =B(u)\text{ with }\lim \sup_{u\rightarrow 0}B(u)<+\infty
\end{equation*}

\noindent We have 
\begin{equation*}
\varphi (v/\sqrt{a_{n}(t_{j})}))-1=-\frac{v^{2}}{2a_{n}(t_{j})}%
+a_{n}^{-3/2}(t_{j})A_{n}=d_{n}\rightarrow 0\text{ with}.\lim
\sup_{n\rightarrow +\infty 0}\left\vert A_{n}\right\vert <+\infty
\end{equation*}

\noindent Then 
\begin{equation*}
\log \varphi (v/\sqrt{a_{n}(t_{j})}))=\log (1+d_{n})=-\frac{v^{2}}{%
2a_{n}(t_{j})}+a_{n}^{-3/2}(t_{j})C_{n}\text{ with }\lim \sup \left\vert
C_{n}\right\vert <+\infty .
\end{equation*}

\noindent Next 
\begin{eqnarray*}
\psi _{Y_{j,n}/\sqrt{a_{n}(t_{j})}}(v) &=&E\exp \left( \left\{ \frac{\Delta
_{n}N(t_{j})}{a_{n}(t_{j})}\right\} \left\{ -\frac{v^{2}}{2}+\frac{C_{n}}{%
a_{n}^{1/2}(t_{j})}\right\} \right) \\
&=&E\exp \left( \left\{ \frac{\Delta _{n}N(t_{j})}{a_{n}(t_{j})}\right\}
\left\{ -\frac{v^{2}}{2}\right\} \right) \exp \left( \left\{ \frac{\Delta
_{n}N(t_{j})}{a_{n}(t_{j})}\right\} \left\{ \frac{C_{n}}{a_{n}^{1/2}(t_{j})}%
\right\} \right)
\end{eqnarray*}

\noindent Let us write $C_{n}=R_{n}(\cos A_{n}+i\sin B_{n})$ with $0\leq
\lim \sup R_{n}<+\infty$. We get 
$$
\psi _{Y_{j,n}/\sqrt{a_{n}(t_{j})}}(v)=E\psi _{Y_{j,n}/\sqrt{a_{n}(t_{j})}}(v)
$$

$$
=\mathbb{E}\exp \left(  \left\{ \frac{\Delta _{n}N(t_{j})}{a_{n}(t_{j})}\right\}
\left\{ -\frac{v^{2}}{2}+\left\{ \frac{R_{n}\cos B_{n}}{a_{n}^{1/2}(t_{j})}
\right\} \right\} \exp \left( i\left\{ \frac{\Delta _{n}N(t_{j})}{a_{n}(t_{j})}\right\} \left\{ \frac{R_{n}\sin A_{n}}{a_{n}^{1/2}(t_{j})}%
\right\} \right) \right) .
$$

\noindent For $n$ large enough, we have $| \frac{R_{n}\cos B_{n}}{a_{n}^{1/2}(t_{j})}| \leq v^{2}/4$ and for 
\begin{equation}
Z_{n}(1)=\exp \left( \left\{ \frac{\Delta _{n}N(t_{j})}{a_{n}(t_{j})}%
\right\} \left\{ -\frac{v^{2}}{2}+\left\{ \frac{R_{n}\cos B_{n}}{%
a_{n}^{1/2}(t_{j})}\right\} \right\} \right) \exp \left( i\left\{ \frac{%
\Delta _{n}N(t_{j})}{a_{n}(t_{j})}\right\} \left\{ \frac{R_{n}\sin A_{n}}{%
a_{n}^{1/2}(t_{j})}\right\} \right) ,  \label{approx0}
\end{equation}

\noindent we have 
\begin{eqnarray}
\left\Vert Z_{n}(1)\right\Vert &=&\exp \left( \left\{ \frac{\Delta
_{n}N(t_{j})}{a_{n}(t_{j})}\right\} \left\{ -\frac{v^{2}}{2}+\left\{ \frac{%
R_{n}\cos B_{n}}{a_{n}^{1/2}(t_{j})}\right\} \right\} \right) \\
&\leq &\exp \left\{ \frac{\Delta _{n}N(t_{j})}{a_{n}(t_{j})}\right\} \left\{
-\frac{v^{2}}{4}\right\} =Z_{n}(2).  \notag
\end{eqnarray}

\noindent Now $Z_{n}$ (2) is of the form 
\begin{equation*}
Z_{n}=g\left( \frac{\Delta _{n}N(t_{j})}{a_{n}(t_{j})}\right)
\end{equation*}

\noindent with 
\begin{equation*}
g(x)=\exp (-\frac{v^{2}}{4}x)1_{(x\geq 0)},
\end{equation*}

\noindent is bounded and continuous on $\mathbb{R}_{+}.$ By $(HN2A)$ 
\begin{equation*}
\frac{\Delta _{n}N(t_{j})}{a_{n}(t_{j})}\rightarrow _{P}1,
\end{equation*}

\noindent which implies that (see Theorem 2.7 in \cite{vaart_asymp}, page
10) 
\begin{equation*}
\frac{\Delta _{n}N(t_{j})}{a_{n}(t_{j})}\rightarrow _{w} 1,
\end{equation*}

\noindent where $\rightarrow _{w}$ stands for weak convergence. By using the
boundedness and the continuity $g$ on $\mathbb{R}_{+}^{2}$ and by the
Portmanteau Theorem or simply by the definition of the weak convergence, we
get 
\begin{equation}
\mathbb{E}Z_{n}(2)=Eg\left( \frac{\Delta _{n}N(t_{j})}{a_{n}(t_{j})}\right)
\rightarrow g(1)=\exp (-\frac{v^{2}}{4})<+\infty .  \label{convA}
\end{equation}

\noindent From (\ref{approx0}), we have 
\begin{equation*}
Z_{n}(1)\rightarrow _{\mathbb{P}}\exp (-\frac{v^{2}}{2}).
\end{equation*}

\noindent We use the Young version of the Dominated Convergence Theorem (see 
\cite{loeve}, page 164) to to get that

\begin{equation*}
\psi _{Y_{j,n}/\sqrt{a_{n}(t_{j})}}(v)=EZ_{n}(1)\rightarrow \exp (-\frac{%
v^{2}}{2}).
\end{equation*}

\noindent Then each $Y_{j,n}/\sqrt{a_{n}(t_{j})}$ weak converges to $N(0,1)$%
. We deduce from this that
\begin{equation*}
(Y_{1,n}/\sqrt{a_{n}(t_{1})},...,Y_{k,n}/\sqrt{a_{n}(t_{k})})\rightarrow
_{d}N_{k}(0,I_{k}),
\end{equation*}

\noindent where $I_{k}$ is identity matrix of dimension $k$. Hence by $(HN3A)
$

\begin{equation*}
(Y_{1,n}/\sqrt{c_{n}},...,Y_{k,n}/\sqrt{c_{n}})\rightarrow
_{d}N_{k}(0,\Lambda _{k}),
\end{equation*}

\noindent where $\Lambda _{k}=diag(a(t_{1}),...,a(t_{k})).$ This leads to %
\ref{fdl01}.\newline

\bigskip \noindent In other words, under the assumptions of this
proposition, the finite distribution of the stochastic process $\{S_{n}(t)/%
\sqrt{c_{n}},0\leq t\leq T\}$ weakly converge to those of a centered
Gaussian process $\{\mathbb{G}(t),0\leq t\leq T\}$ with covariance function

\begin{equation*}
\Gamma (s,t)=a(\min (s,t)),
\end{equation*}

\noindent where

\begin{equation*}
a(s,t)=\lim_{n\rightarrow +\infty }\mathbb{E}(N(nt)-N(ns))/c_{n}.
\end{equation*}

\bigskip \noindent \textbf{Finally}, let us address the tightness of the
sequence.

\begin{lemma}
\label{lem02} Suppose that the assumptions of Proposition \ref{prop01} hold.
Assume in addition that $(HTA)$ holds. Then the stochastic process 
$$
\{S_{n}(t)/\sqrt{c_{n}},0\leq t\leq T\}
$$ 

\noindent is asymptotically tight in $\ell
^{\infty }(T)$\bigskip , that is (\ref{formulaTensIndep}) holds.
\end{lemma}

\bigskip \noindent \textbf{Proof of Lemma \ref{lem02}}. Fix $\eta >0$ and $%
t\in \lbrack 0,T[$ and take $\delta >0$ such that $t+\delta \leq T.$ For $%
s\in [t,t+\delta]$, 
\begin{equation*}
\frac{\left\vert S_{n}(s)-S_{n}(t)\right\vert }{\sqrt{c_{n}}}=\frac{1}{\sqrt{%
c_{n}}}\left\vert \sum_{h=N(nt)+1}^{N(ns)}X_{h}(\omega )\right\vert .
\end{equation*}

\noindent We may and do replace all the process $\left\{ Y_{n}(s),s\in
[t,t+\delta] \right\} $ by 
\begin{equation*}
Y_{n}(s)=\left\{ \frac{1}{\sqrt{c_{n}}}\left\vert
\sum_{h=1}^{N(ns)-N(nt)}X_{h}\right\vert ,s\in [t,t+\delta] \right\}
\end{equation*}

\noindent in the sense of equality in law since the random variables $%
X_{1},X_{2},...$ are iid. Next, since that random variable are integers, the
supremum 
\begin{equation*}
Y_{n}(t,\delta )=\sup_{s\in \lbrack [t,t+\delta]}Y_{n}(s)
\end{equation*}

\noindent is taken over a countable number of values, and then is
measurable. Denoting $\Delta _{n}N(t,\delta )=N(nt+n\delta )-N(nt),$ we see
that%
\begin{equation*}
Y_{n}(t,\delta )\leq \frac{1}{\sqrt{c_{n}}}\max_{j\leq \Delta _{n}N(t,\delta
)}\left\vert \sum_{h=1}^{j}X_{h}\right\vert .
\end{equation*}

\noindent To simplify, denote $S_{0}=0,$ $S_{j}=\sum_{h=1}^{j}X_{h}(\omega
),j\geq 1.$ we have 
\begin{equation*}
Y_{n}(t,\delta )\leq \frac{1}{\sqrt{c_{n}}}\max_{j\leq \Delta _{n}N(t,\delta
)}\left\vert S_{j}\right\vert .
\end{equation*}

\noindent From now, we do not need the use of outer probability since $%
Y_{n}(t,\delta )$ is measurable. We have

\begin{eqnarray}
{\large \mathbb{P}(}Y_{n}(t,\delta ) &>&{\large \eta )\leq \mathbb{P}(}%
\max_{j\leq \Delta _{n}N(t,\delta )}\left\vert S_{j}\right\vert >{\large %
\eta \sqrt{c_{n}})}  \label{etapeM4} \\
&\leq &\sum_{r=0}^{\infty }{\large \mathbb{P}(}\Delta _{n}N(t,\delta )=r%
{\large )\mathbb{P}(}\max_{j\leq r}\left\vert S_{j}\right\vert >{\large \eta 
\sqrt{c_{n}}).}  \notag \\
&=&\sum_{r=1}^{\infty }{\large \mathbb{P}(}\Delta _{n}N(t,\delta )=r{\large )%
\mathbb{P}(}\max_{j\leq r}\left\vert S_{j}\right\vert >{\large \eta \sqrt{%
c_{n}})}.  \notag
\end{eqnarray}

\noindent For any fixed $r\geq 1$, for $p\geq 1,$ the sequence $\left\vert
S_{1}\right\vert ^{p},...,\left\vert S_{r}\right\vert ^{p}$ is a
submartingale and then satisfies the maximal inequality 
\begin{equation*}
{\large \mathbb{P}(}\max_{j\leq r}\left\vert S_{j}\right\vert >{\large \eta
)=\mathbb{P}(\max_{j\leq r}\left\vert S_{j}\right\vert }^{p}{\large >\eta }%
^{p}{\large )\leq }\frac{E\left\vert S_{r}\right\vert ^{p}}{c_{n}^{p/2}\eta
^{p}}.
\end{equation*}

\noindent Then 
\begin{eqnarray*}
{\large \mathbb{P}(}Y_{n}(t,\delta ) &>&{\large \eta )\leq }%
\sum_{r=1}^{\infty }{\large \mathbb{P}(}\Delta _{n}N(t,\delta )=r{\large )}%
\frac{E\left\vert S_{r}\right\vert ^{p}}{c_{n}^{p/2}\eta ^{p}} \\
&=&\frac{c_{n}^{-p/2}}{\eta ^{p}}\sum_{r=1}^{\infty }{\large \mathbb{P}(}%
\Delta _{n}N(t,\delta )=r{\large )r}^{p/2}\left( E\left\vert \frac{S_{r}}{%
\sqrt{r}}\right\vert ^{p}\right) .
\end{eqnarray*}

\noindent Since the sequence 
\begin{equation*}
\mathbb{E}\left( \left\vert \frac{S_{r}}{\sqrt{r}}\right\vert ^{p}\right)
,p\geq 1,
\end{equation*}

\noindent is bounded, say by $C>0,$ we get 
\begin{eqnarray*}
{\large \mathbb{P}(}Y_{n}(t,\delta ) &>&{\large \eta )\leq }\frac{%
Cc_{n}^{-p/2}}{\eta ^{p}}\sum_{r=1}^{\infty }{\large \mathbb{P}(}\Delta
_{n}N(t,\delta )=r{\large )r}^{p/2} \\
&=&\frac{Cc_{n}^{-p/2}}{\eta ^{p}}\mathbb{E}{\large (}\Delta _{n}N(t,\delta
)^{p/2}).
\end{eqnarray*}

\noindent Then 
\begin{equation*}
\lim_{\delta \rightarrow 0}\sup_{s\in \lbrack 0,T]}\limsup_{n}\frac{1}{%
\delta } \mathbb{P}(\sup_{s-\delta <t<s+\delta ,t\in \lbrack 0,T]}\frac{%
\left\vert S_{n}(s)-S_{n}(t)\right\vert }{\sqrt{c_{n}}}>\eta) 
\end{equation*}

\begin{equation*}
\leq \lim_{\delta \rightarrow 0}\lim \sup_{n\rightarrow +\infty }\frac{C}{%
\eta ^{p}\delta }\mathbb{E}\left( \left\{ \frac{\Delta _{n}N(t,\delta )}{%
c_{n}}\right\} ^{p/2}\right). 
\end{equation*}

\noindent Then the sequence is asymptotically tight whenever $(HTA)$ holds
for some $p\geq 1$.\newline

\noindent We conclude the proof of Proposition \ref{propA} by combining
Proposition \ref{prop01} and Lemma \ref{lem02}, we get the searched result.\newline

\bigskip 
\subsection{Proof of Corollary ref{corA}}

Suppose that the application $a(\circ )$ transforms $[0,T]$ into $[0,A_{T}]$
with $a(0)=0.$ Suppose $a(\circ )$ is invertible, and $a^{-1}(\circ )$ is $k$%
-Lipschitz. We may apply Proposition \ref{prop01} to $\left\{
S_{n,T}(a^{-1}(u),0\leq u\leq A_{T}\right\} .$ For 
\begin{equation*}
0=t_{0}=a(u_{0})<t_{2}=a(u_{2})<...<t_{k}=a(u_{k}),
\end{equation*}

\noindent we surely have that
\begin{equation*}
\left( \frac{S_{n,T}(a^{-1}(u_{1})-S_{n,T}(a^{-1}(u_{0})}{\sqrt{c_{n}}},...,%
\frac{S_{n,T}(a^{-1}(u_{k})-S_{n,T}(a^{-1}(u_{k-1})}{\sqrt{c_{n}}}\right) 
\end{equation*}

\noindent weakly converges to ($u_{0}=0$,

\begin{equation} 
(\Delta B(1)\sqrt{a(a(a^{-1}(u_{1}))}+...+\Delta B(j)\sqrt{a(a(a^{-1}(u_{j}))},1\leq j \leq k),  \label{FDT}
\end{equation}

\noindent where $a(a^{-1}(u_{j}))=a(a^{-1}(u_{j-1})$, $a^{-1}(u_{j}))=a(a^{-1}(u_{j}))-a(a^{-1}(u_{j-1}))=u_{j}-u_{j-1}$. It comes that (\ref{FDT}), which is
\begin{equation*}
( \Delta B(1)\sqrt{u_{1}}+\Delta B(2)\sqrt{u_{2}-u_{1}}+...+\Delta B(j)\sqrt{u_{k}-u_{k-1}},1\leq j \leq k),  \label{FDT}
\end{equation*}

\noindent is a finite distribution of a the Brownian motion $\{B,0\leq u\leq A_{T}\}.$\\

\bigskip \noindent Its remains to check that $S_{n,T}(a^{-1}(\circ ))$ is asymptotically tight under the assumptions. We have for any fixed $v=a(t)\in ]0,A_{T}[$,  for $\delta>0$ such that $v+\delta \leq A_{T}$ and $a^{-1}(v)+k\delta \leq T$
for any $u$ such that $v<u<v+\delta $ 
\begin{equation*}
\frac{\left\vert S_{n}(a^{-1}(u))-S_{n}(a^{-1}(v))\right\vert }{\sqrt{c_{n}}}%
=\frac{\left\vert S_{n}(s)-S_{n}(t)\right\vert }{\sqrt{c_{n}}},
\end{equation*}

\noindent where $s=a^{-1}(u)$ and $s-t=a^{-1}(u)-a^{-1}(v)\leq k(u-v)\leq k\delta $
and then

\begin{equation*}
\frac{\left\vert S_{n}(a^{-1}(u))-S_{n}(a^{-1}(v))\right\vert }{\sqrt{c_{n}}}%
\leq \sup_{t\leq s\leq t+k\delta }\frac{\left\vert
S_{n}(s)-S_{n}(t)\right\vert }{\sqrt{c_{n}}},
\end{equation*}

\noindent and next

\begin{equation*}
\sup_{v\leq u\leq v+\delta }\frac{\left\vert
S_{n}(a^{-1}(u))-S_{n}(a^{-1}(v))\right\vert }{\sqrt{c_{n}}}\leq \sup_{t\leq
s\leq t+k\delta }\frac{\left\vert S_{n}(s)-S_{n}(t)\right\vert }{\sqrt{c_{n}}%
},
\end{equation*}

\noindent and next for any $\eta>0$
\begin{equation*}
\mathbb{P}\left( \sup_{v\leq u\leq v+\delta }\frac{\left\vert
S_{n}(a^{-1}(u))-S_{n}(a^{-1}(v))\right\vert }{\sqrt{c_{n}}}>\eta \right)
\leq \mathbb{P}\left( \sup_{t\leq s\leq t+k\delta }\frac{\left\vert
S_{n}(s)-S_{n}(t)\right\vert }{\sqrt{c_{n}}}>\eta \right) .
\end{equation*}

\noindent We conclude that for any fixed $v\in ]0,A_{T}[.$ 
\begin{equation*}
\lim_{\delta \rightarrow 0}\frac{1}{\delta }\lim \sup_{n\rightarrow +\infty }%
\mathbb{P}\left( \sup_{v\leq u\leq v+\delta }\frac{\left\vert
S_{n}(a^{-1}(u))-S_{n}(a^{-1}(v))\right\vert }{\sqrt{c_{n}}}>\eta \right) =0.
\end{equation*}

\noindent The sequence $S_{n,T}(a^{-1}(\circ ))$ is tight. Then, in combination of the convergence of its the finite distributions to those of a Brownian motiona, it weakly converges to the Brownian motion on  $[0,A_{T}]$


\begin{thebibliography}{9}
\bibitem{loeve} Lo\`{e}ve, Michel.(1997). Probability measure.
Springer-Verlag, 4th Edition.

\bibitem{vaart} A. W. van der Vaart and J. A. Wellner(1996). \textit{Weak
Convergence and Empirical Processes With Applications to Statistics}.
Springer, New-York.

\bibitem{vaart_asymp} van der Vaart, A.W. \textit{Asymptotic Statistics}.
(2000). \textit{Cambridge}

\bibitem{billingsley} Billingsley, P. (1968). \textit{Convergence of
probability measures}. Wiley, New-York.

\bibitem{loTensTool} Lo, G.S. (2014). A remark on the asymptotic tightness
in $\ell^{\infty}([a,b])$. arxiv.org/pdf/1405.6342.
\end{thebibliography}
\end{document}